\documentclass[11pt,fleqn,twoside]{article}
\usepackage{amsfonts,amssymb,latexsym}
\makeatletter
\newcommand{\prava}[1]{\small\it
\begin{flushleft}
Copyright \copyright \ 1999 by  #1
\end{flushleft}}

\newcommand{\name}[1]{\begin{flushleft}
                       \LARGE \bf #1
                       \end{flushleft}\vspace{-3mm}}

\newcommand{\Author}[1]{\begin{flushleft}
                       \it #1 \end{flushleft}}

\newcommand{\Adress}[1]{\begin{flushleft}
                       \it #1 \end{flushleft}}

\newcommand{\Date}[1]{\begin{flushleft}
                      \small  \it #1 \end{flushleft}}

\newcommand{\ehkol}{Author \ name}
\newcommand{\ohkol}{Article \ name}
\renewcommand{\@evenhead}{
\hspace*{-3pt}\raisebox{-15pt}[\headheight][0pt]{\vbox{\hbox to \textwidth 
{\thepage \hfil \ehkol}\vskip4pt \hrule}}}
\renewcommand{\@oddhead}{
\hspace*{-3pt}\raisebox{-15pt}[\headheight][0pt]{\vbox{\hbox to \textwidth 
{\ohkol \hfil \thepage}\vskip4pt\hrule}}}
\renewcommand{\@evenfoot}{}
\renewcommand{\@oddfoot}{}

\newcommand{\be}{\begin{equation}}
\newcommand{\ee}{\end{equation}}
\newcommand{\ba}{\hspace*{-5pt}\begin{array}}
\newcommand{\ea}{\end{array}}

\newcommand{\ds}{\displaystyle}
\makeatother

\begin{document}
\thispagestyle{empty}
\setcounter{page}{120}
\renewcommand{\ehkol}{R. Radha, S. Vijayalakshmi and M. Lakshmanan}
\renewcommand{\ohkol}{Explode-Decay Dromions in the Non-Isospectral
Davey-Stewartson I Equation}

\begin{flushleft}
\footnotesize \sf
Journal of Nonlinear Mathematical Physics \qquad 1999, V.6, N~2,
\pageref{radha-fp}--\pageref{radha-lp}.
\hfill {\sc Letter}
\end{flushleft}

\vspace{-5mm}

\renewcommand{\footnoterule}{}
{\renewcommand{\thefootnote}{}
\footnote{\prava{R. Radha, S. Vijayalakshmi and M. Lakshmanan}}}

\name{Explode-Decay Dromions in the Non-Isospectral Davey-Stewartson I
(DSI) Equation}\label{radha-fp}

\Author{R. RADHA~$^{\dag}$, S. VIJAYALAKSHMI~$^{\ddag}$ and M. LAKSHMANAN~$^{\ddag}$}

\Adress{$\dag$~Department of Physics, Government College for Women,\\
~~Kumbakonam -- 612 001, India\\[1mm]
$\ddag$~Centre for Nonlinear Dynamics, Department of Physics,\\
~~Bharathidasan University, Tiruchirapalli 620 024, India\\
~~E-mail: lakshman@bdu.ernet.in}

\Date{Received October 10, 1998; Revised February 25, 1999;
Accepted March 25, 1999}

\begin{abstract}
\noindent
In this letter, we report the existence of a novel type of
explode-decay dromions, which are exponentially localized coherent
structures whose amplitude varies with time, through Hirota method for
a nonisospectral Davey-Stewartson equation~I discussed recently by
Jiang. Using suitable transformations, we also point out such solutions
also exist for the isospectral Davey-Stewartson~I equation itself for a
careful choice of the potentials.
\end{abstract}

The identif\/ication of exponentially localized structures namely the so called
``dromions" in the Davey-Stewartson~I (DSI) equation~[1] has given a new direction
to the study of nonlinear partial dif\/ferential equations (pdes) in (2+1)
dimensions. The subsequent investigations of other physically and mathematically
important (2+1) dimensional nonlinear pdes [2--10] have enriched the family of
integrable nonlinear pdes, exhibiting considerable richness in the structure
of the solutions. However, very little is known about non-isospectral (2+1)
dimensional nonlinear pdes in this regard. Now one may ask, do exponentially
localized structures exist in the later systems also, and if so,
what are their
characteristics? This letter is a modest attempt in this direction in which the
existence of such localized solutions is pointed out for DSI equation with
nonuniform or inhomogeneous terms.

In this letter, we consider a non-isospectral DSI equation introduced by Jiang~[11] of the form
\renewcommand{\theequation}{\arabic{equation}{\rm a}}
\be
\ba{l}
\ds iq_t+{1\over 2}(q_{xx}+q_{yy})+\left[ u_1(\xi,t)+ u_2(\eta,t)-{1\over 2}
\partial_{\xi}^{-1}(qr)_{\eta}-{1 \over 2}\partial_{\eta }^{-1}(qr)_{\xi }\right]q
\vspace{3mm}\\
\ds \qquad - i(\omega_1y+a_1)q_y-i\omega_1(xq)_x+2(\omega_0x+a_0)q = 0,
\ea
\ee
\renewcommand{\theequation}{\arabic{equation}{\rm b}}
\setcounter{equation}{0}
\be
\ba{l}
\ds ir_t-{1\over 2}(r_{xx}+r_{yy})-\left[ u_1(\xi,t)+ u_2(\eta,t)-{1\over 2}
\partial_{\xi}^{-1}(qr)_{\eta}-\partial_{\eta }^{-1}(qr)_{\xi }\right]r
\vspace{3mm}\\
\ds \qquad - i(\omega_1y+a_1)r_y-i\omega_1(xr)_x+
2(\omega_0x+a_0)r = 0,\qquad  \xi = x+y,\quad \eta = x-y,
\ea\hspace{-4.8pt}
\ee
where $\omega_0$, $\omega_1$, $a_0$ and $a_1$ are all real
constants and $u_1(\xi,t)$ and $u_2(\eta,t)$ are the so called boundaries
(arbitrary functions). The above set of equations has been investigated
via inverse scattering transform method~[11] and an integro-dif\/ferential
equation for the time evolution of the scattering data  by virtue of the
time dependence of the scattering parameters has been brought out. Under
the reduction $r = -q^*$, equation~(1) reduces to the following form
\renewcommand{\theequation}{\arabic{equation}}
\setcounter{equation}{1}
\be
\ba{l}
\ds iq_t+q_{\xi \xi}+q_{\eta \eta}+\left[ u_1(\xi,t)+ u_2(\eta,t)+{1\over 2}
\partial_{\xi}^{-1}(\mid q\mid ^2)_{\eta}+{1\over 2}\partial_{\eta }^{-1}(\mid q\mid ^2)_
{\xi}\right] q
\vspace{3mm}\\
\ds \qquad -i\omega_1[\xi q_{\xi}+\eta q_{\eta}]-ia_1[q_{\xi}-q_{\eta}]
+[2a_0+\omega_0(\xi+\eta)-i\omega_1]q = 0.
\ea
\ee
In equation (2), the parameter $a_0$ can be scaled away by
introducing the transformation $q = \hat q \exp {(-2ia_0t)}$. Then, the above
equation can be equivalently written as
\renewcommand{\theequation}{\arabic{equation}{\rm a}}
\setcounter{equation}{2}
\be
{1\over 2}\left(\mid q\mid ^2\right)_{\eta} = U_{\xi},
\ee
\renewcommand{\theequation}{\arabic{equation}{\rm b}}
\setcounter{equation}{2}
\be
{1\over 2}\left(\mid q\mid
^2\right)_{\xi} = V_{\eta},
\ee
\renewcommand{\theequation}{\arabic{equation}}
\setcounter{equation}{3}
\be
iq_t+q_{\xi\xi}+q_{\eta\eta}+(U+V)q-i\omega_1(\xi q_\xi+\eta q_\eta)
-ia_1(q_\xi-q_\eta)+(\omega_0(\xi+\eta)-i\omega_1)q = 0.
\ee
As the complete integrability parallelling to that of dynamical systems
under isospectral f\/lows is not obvious for the nonisospectral problems
in (2+1) dimensions, we address ourselves only to the nature of the
solutions of the above equation rather than its integrability property. For
this purpose, we bilinearize equation~(3) and obtain the solutions using the
Hirota method.

To bilinearise the above equation, we ef\/fect the following dependent
variable transformation
\be
q = {G\over F},\qquad U = 2\partial_{\eta\eta}\log F,\qquad
V = 2\partial_{\xi\xi}\log F,
\ee
so that equations (3) and (4) get converted into the following Hirota form,
\renewcommand{\theequation}{\arabic{equation}{\rm a}}
\setcounter{equation}{5}
\be
\ba{l}
\left[ iD_t+D_{\xi}^2+D_{\eta}^2-i\omega_1(\xi D_{\xi}+\eta D_{\eta})\right.
\vspace{2mm}\\
\ds \qquad \left. -ia_1(D_{\xi}-D_{\eta})+(\omega_0[\xi+\eta]-i\omega_1) \right]G \cdot F = 0,
\ea
\ee
\renewcommand{\theequation}{\arabic{equation}{\rm b}}
\setcounter{equation}{5}
\be
2D_\xi D_\eta F\cdot F = \mid G\mid ^2.
\ee

We now introduce the following power series expansion
\renewcommand{\theequation}{\arabic{equation}{\rm a}}
\setcounter{equation}{6}
\be
G = \epsilon g^{(1)}+\epsilon^3 g^{(3)}+ \cdots,
\ee
\renewcommand{\theequation}{\arabic{equation}{\rm b}}
\setcounter{equation}{6}
\be
F = 1+\epsilon^2 f^{(2)}+\epsilon^4 f^{(4)}+ \cdots,
\ee
into the bilinear form (6), where $\epsilon$ is a small parameter. Collecting
the various powers of~$\epsilon$, we get the following set of equations,
\renewcommand{\theequation}{\arabic{equation}{\rm a}}
\setcounter{equation}{7}
\be
\ba{ll}
O(\epsilon): & \ i(g^{(1)})_t+(g^{(1)})_{\xi\xi}+(g^{(1)})_{\eta\eta}-i\omega_1[\xi(g^{(1)})_{\xi}+
\eta(g^{(1)})_{\eta}]
\vspace{2mm}\\
&   -ia_1[(g^{(1)})_{\xi} -(g^{(1)})_{\eta}] +(\omega_0(\xi+\eta)-i\omega_1) g^{(1)} = 0,
\ea
\ee
\renewcommand{\theequation}{\arabic{equation}{\rm b}}
\setcounter{equation}{7}
\be
\ba{ll}
O(\epsilon^2): & 4(f^{(2)})_{\xi\eta} = g^{(1)} g^{(1)*},
\ea
\ee
\renewcommand{\theequation}{\arabic{equation}{\rm c}}
\setcounter{equation}{7}
\be
\ba{ll}
O(\epsilon^3): & \biggl [ iD_t+D_{\xi}^2+D_{\eta}^2-i\omega_1(\xi D_{\xi}+\eta D_{\eta})
\vspace{2mm}\\
 &  -ia_1(D_{\xi}-D_{\eta})+(\omega_0[\xi+\eta]-i\omega_1) \biggr]+
\left(g^{(3)}+g^{(1)}.f^{(2)}\right) = 0,
\ea
\ee
\renewcommand{\theequation}{\arabic{equation}{\rm d}}
\setcounter{equation}{7}
\be
\ba{ll}
O(\epsilon^4): & 2D_\xi D_\eta(2f^{(4)}+f^{(2)}.f^{(2)}) = g^{(3)}g^{(1)^*}+g^{(1)}g^{(3)^*},
\ea
\ee
and so on. Solving (8a), we obtain the simplest ``plane wave'' solution
\renewcommand{\theequation}{\arabic{equation}}
\setcounter{equation}{8}
\be
g^{(1)} = \sum_{j = 1}^N e^{\chi_j},\qquad
\chi_j = k_j(t)\xi+l_j(t)\eta+\int\Omega_j(t)dt,
\ee
where the spectral parameters $k_j(t)$ and $l_j(t)$ evolve in an identical
fashion obeying the time evolution equation
\be
i(k_j)_t-i(\omega_1)k_j = i(l_j)_t-i(\omega_1)l_j = -\omega_0,
\ee
and the ``dispersion'' relation is given by
\be
i\Omega_j(t)+(k_j(t))^2+(l_j(t))^2-ia_1(k_j(t)-l_j(t))-i\omega_1 = 0.
\ee

To construct one soliton solution, we take $N = 1$ so that we have
\be
g^{(1)} = e^{\chi_1},
\ee
and hence the solution of (8b) becomes
\be
f^{(2)} = e^{\chi_1+\chi_1^*+2\psi},\qquad
e^{2\psi} = {1\over 16k_{1R}(t)l_{1R}(t)}.
\ee
Substituting $g^{(1)}$ and $f^{(2)}$ in equations (8c) and (8d), one can show that
$g^{(j)} = 0$ for $j\ge3$ and $f^{(j)} = 0$ for $j\ge4$ provided the spectral
parameters obey the following equation
\be
(k_{1R}l_{1R})_t = 2\omega_1 k_{1R}l_{1R}.
\ee
Considering the above equation alongwith equation (10), separating out the real
and imaginary parts of the later, the time evolution of the spectral
parameters can be obtained as
\renewcommand{\theequation}{\arabic{equation}{\rm a}}
\setcounter{equation}{14}
\be
k_{1R}(t) = k_{1R}(0)e^{\omega_1t},\qquad l_{1R}(t) = l_{1R}(0)e^{\omega_1t}
\ee
and
\renewcommand{\theequation}{\arabic{equation}{\rm b}}
\setcounter{equation}{14}
\be
k_{1I}(t) = k_{1I}(0)e^{\omega_1t}-\frac{\omega_0}{\omega_1},\qquad
l_{1I}(t) = l_{1I}(0)e^{\omega_1t}-\frac{\omega_0}{\omega_1},
\ee
where $k_{1R}(0)$, $k_{1I}(0)$, $l_{1R}(0)$ and $l_{1I}(0)$ are all
constants. Using (5), (7), (12) and (13), the physical f\/ield variable $q$
of equation~(4) is driven by the envelope soliton (line soliton)
\renewcommand{\theequation}{\arabic{equation}}
\setcounter{equation}{15}
\be
q =2\sqrt{k_{1R}(t) l_{1R}(t)} \, \mbox{sech}\,  (\chi_{1R}+\psi)e^{i\chi_{1I}}.
\ee
It is interesting to note that the amplitude of the above soliton
solution varies with time by virtue of equation~(15). Similarly, the
potentials $U$ and $V$ are driven by the line soliton solutions whose
amplitude again varies with time as
\renewcommand{\theequation}{\arabic{equation}{\rm a}}
\setcounter{equation}{16}
\be
U = (k_{1R}(t))^2 \, \mbox{sech}^2(\chi_{1R}+\psi),
\ee
\renewcommand{\theequation}{\arabic{equation}{\rm b}}
\setcounter{equation}{16}
\be
V = (l_{1R}(t))^2 \, \mbox{sech}^2(\chi_{1R}+\psi).
\ee

As it is clear from equation (3) that the boundaries are specif\/ied by the
arbitrary functions $u_2(\eta,t$) and $u_1(\xi,t)$ which drive the potentials
$U$ and $V$ even in the absence of the f\/ield $q$, as in the case of isospectral
(uniform) DSI equation, one can expect the presence of exponentially localized
solutions in the system provided one properly takes care of the time evolution
of the spectral parameters. To generate a (1,1) dromion solution, we take the ansatz
\renewcommand{\theequation}{\arabic{equation}}
\setcounter{equation}{17}
\be
\ba{l}
\ds F = \delta+\alpha e^{\chi_1+\chi_1^*}+\beta e^{\chi_2+\chi_2^*}+
\gamma e^{\chi_1+\chi_1^*+\chi_2+\chi_2^*}, \qquad
  \chi_1 = k_1\xi+ \int\Omega_1(t) dt,
\vspace{3mm}\\
\ds \chi_2 = l_1\eta+\int\Omega_2(t)dt, \qquad
\Omega_1(t) = ik_1^2+a_1k_1,\qquad \Omega_2(t) = il_1^2-a_1l_1.
\ea
\ee
where $\alpha$, $\beta$, $\gamma$ and $\delta$ are parameters.
Substituting (18) into (6b), we obtain
\renewcommand{\theequation}{\arabic{equation}{\rm a}}
\setcounter{equation}{18}
\be
G = \rho e^{\chi_1+\chi_2},\qquad \rho \quad \mbox{is real},
\ee
\renewcommand{\theequation}{\arabic{equation}{\rm b}}
\setcounter{equation}{18}
\be
\rho^2 = 8k_{1R}(t) l_{1R}(t) (\delta\gamma-\alpha\beta).
\ee
Hence, the dromion solution now becomes
\renewcommand{\theequation}{\arabic{equation}}
\setcounter{equation}{19}
\be
q = {{\rho e^{\chi_1+\chi_2}}\over {\delta+\alpha e^{\chi_1+\chi_1^*}
+\beta e^{\chi_2+\chi_2^*}+\gamma e^{\chi_1+\chi_1^*+\chi_2+\chi_2^*}}}.
\ee
It can be easily observed from the above solution that the amplitude of
the dromion solution has to evolve in time obeying the time evolution equation
\be
\rho_t-\omega_1\rho = 0,
\ee
so that solution (20) satisf\/ies equation~(6). The above equation has the solution
\be
\rho = \rho_0 e^{\omega_1t},
\ee
where $\rho_0$ is a constant. Thus, the amplitude of the dromion
solution varies with time in the nonisospectral DSI equation unlike the
isospectral DSI equation where it remains a constant. It should be
mentioned that both the spectral parameters $k_j$ and $l_j$ as well as
the amplitude of the dromion solution are governed by the same time
evolution equation (dif\/fering only in integration constants). Thus,
depending on the nature of the parameters involved, the amplitude of
the dromion solution (20) either grows (explodes) or decays with time.
A typical example is shown in  Figs.~1--3. We call these types of
solutions as ``explode-decay dromions'' which is reminiscent of the
explode-decay solitons of the inhomogeneous (1+1) dimensional nonlinear
pdes~[12].  To our knowledge this seems to be the f\/irst instance of
such a localized solution being realized in a (2+1) dimensional
nonlinear pde.

Finally, we wish to point out that the above type of localized solutions with
time varying amplitude do exist for the isospectral Davey-Stewartson I equation
also. Using the following transformations
\renewcommand{\theequation}{\arabic{equation}{\rm a}}
\setcounter{equation}{22}
\be
\hat q = q \exp {\left(-i\left[{1\over 4} \omega_1\left(\xi^2+\eta^2\right)+{1\over 4}a_1(\xi-\eta)
+{1\over 2} a_1^2t\right]\right)},
\ee
\renewcommand{\theequation}{\arabic{equation}{\rm b}}
\setcounter{equation}{22}
\be
\hat U = U+{1\over 4}\omega_1^2\eta^2-\left({1\over 2}a_1\omega_1-\omega_0\right) \eta,
\ee
\renewcommand{\theequation}{\arabic{equation}{\rm c}}
\setcounter{equation}{22}
\be
\hat V = V+{1\over 4}\omega_1^2\xi^2-\left({1\over 2}a_1\omega_1-\omega_0\right) \xi,
\ee
the nonisospectral DSI equation (3)--(4) reduces to the isospectral DSI equation
\renewcommand{\theequation}{\arabic{equation}{\rm a}}
\setcounter{equation}{23}
\be
i\hat q_t+\hat q_{\xi \xi}+\hat q_{\eta \eta}+\left(\hat U+\hat V\right)\hat q = 0,
\ee
\renewcommand{\theequation}{\arabic{equation}{\rm b}}
\setcounter{equation}{23}
\be
\hat U_{\xi} = {1\over 2} \mid \hat q\mid^2_{\eta},
\ee
\renewcommand{\theequation}{\arabic{equation}{\rm c}}
\setcounter{equation}{23}
\be
\hat V_{\xi} = {1\over 2} \mid \hat q\mid^2_{\xi}.
\ee
It is well known that equation (24) admits exponentially localized solutions
with constant amplitude by driving $\hat U$ and $\hat V$ by $\mbox{sech}^2$ potentials
[1, 8]. However, by virtue of the above transformation (23), it is now evident
that such localized solutions with time varying amplitude do exist for the
isospectral DSI equation also for a careful choice of the potentials indicated
by the transformations (23b) and (23c). However, this choice of the potential
is not very obvious but for the above nonisospectral case by virtue of the
transformations~(23) and hence such localized solutions whose amplitude either
grows or decays with time have eluded earlier observation.

In this letter, we have generated a new class of localized coherent
structures to DSI equation with inhomogeneous terms known as ``explode-decay
dromions" whose amplitude varies with time unlike the basic dromions. It remains
to be seen how the multiexplode-decay dromions would interact in the context of
the variation of the spectral parameter and this remains as an open question. We
have also indicated the possibility of the existence of such solutions for the
isospectral DSI equation itself.

The authors would like to acknowledge the f\/inancial support from the Department
of Science and Technology in the form of a research project. S.V. wishes to
thank Council of Scientif\/ic and Industrial Research, India, for providing a
Senior Research Fellowship.

%\label{radha-lp}

\strut\hfill

{\bf Figure Captions}

%{\bf Figure Captions}

%\vskip2cm

\strut\hfill

\centerline{ {\bf Fig. 1:} The time evolution of explode-decay
dromion at $t = -0.5$}

\vfill

\pagebreak

\vskip2cm
\centerline{ {\bf Fig. 2:} $\mid q \mid$ at $t = 0$}

\vskip12cm

\centerline{ {\bf Fig. 3:} $\mid q \mid$ at $t = 0.2$}

\label{radha-lp}

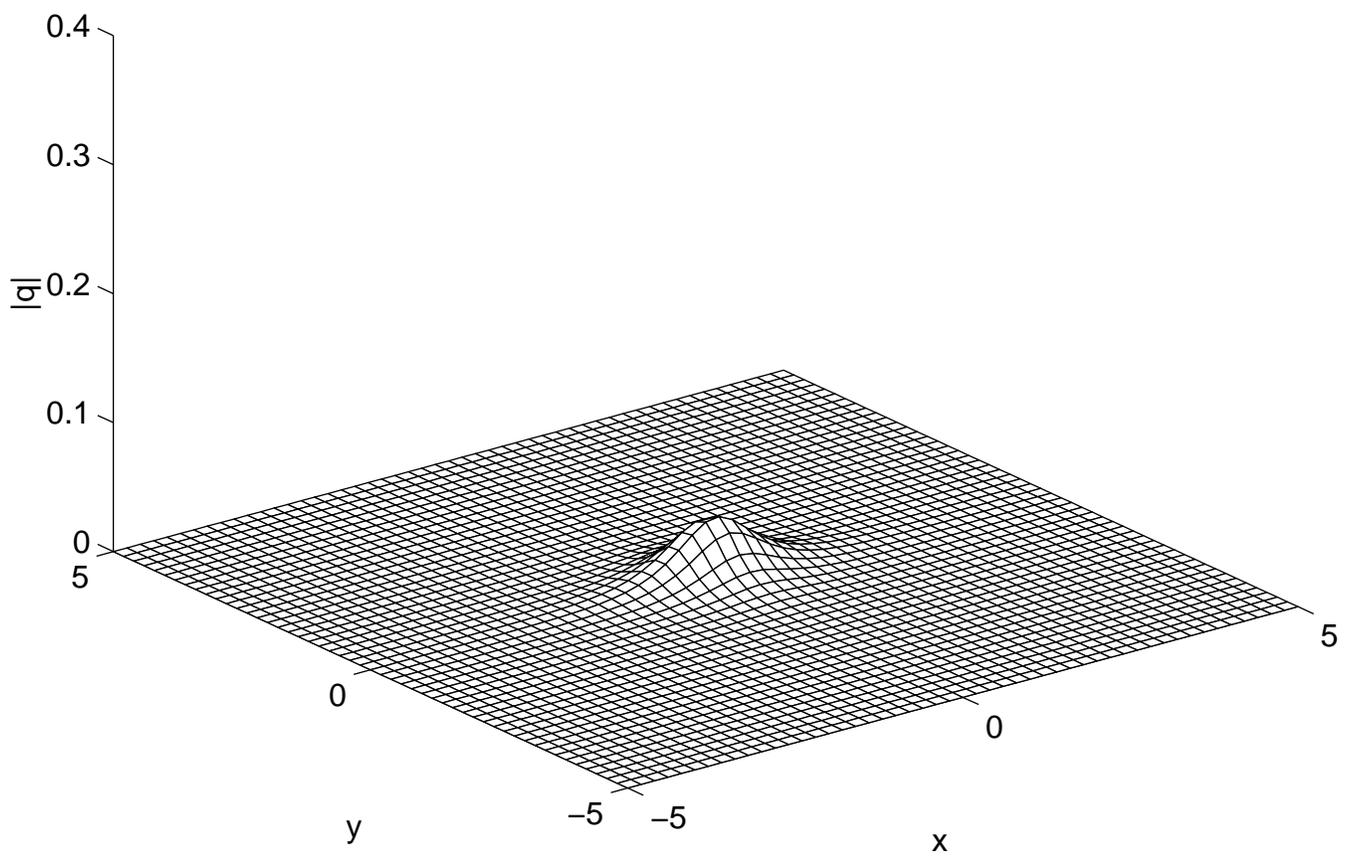
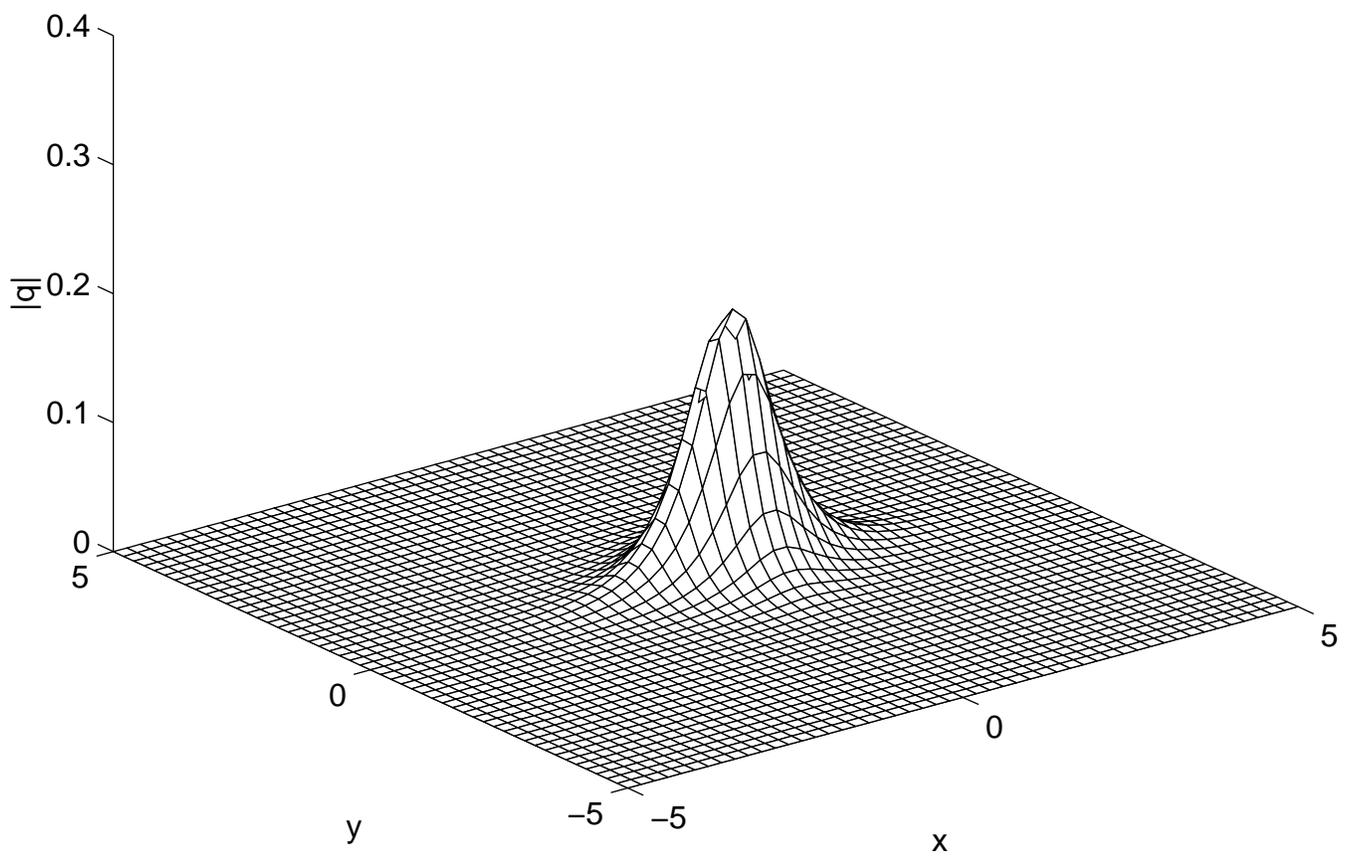
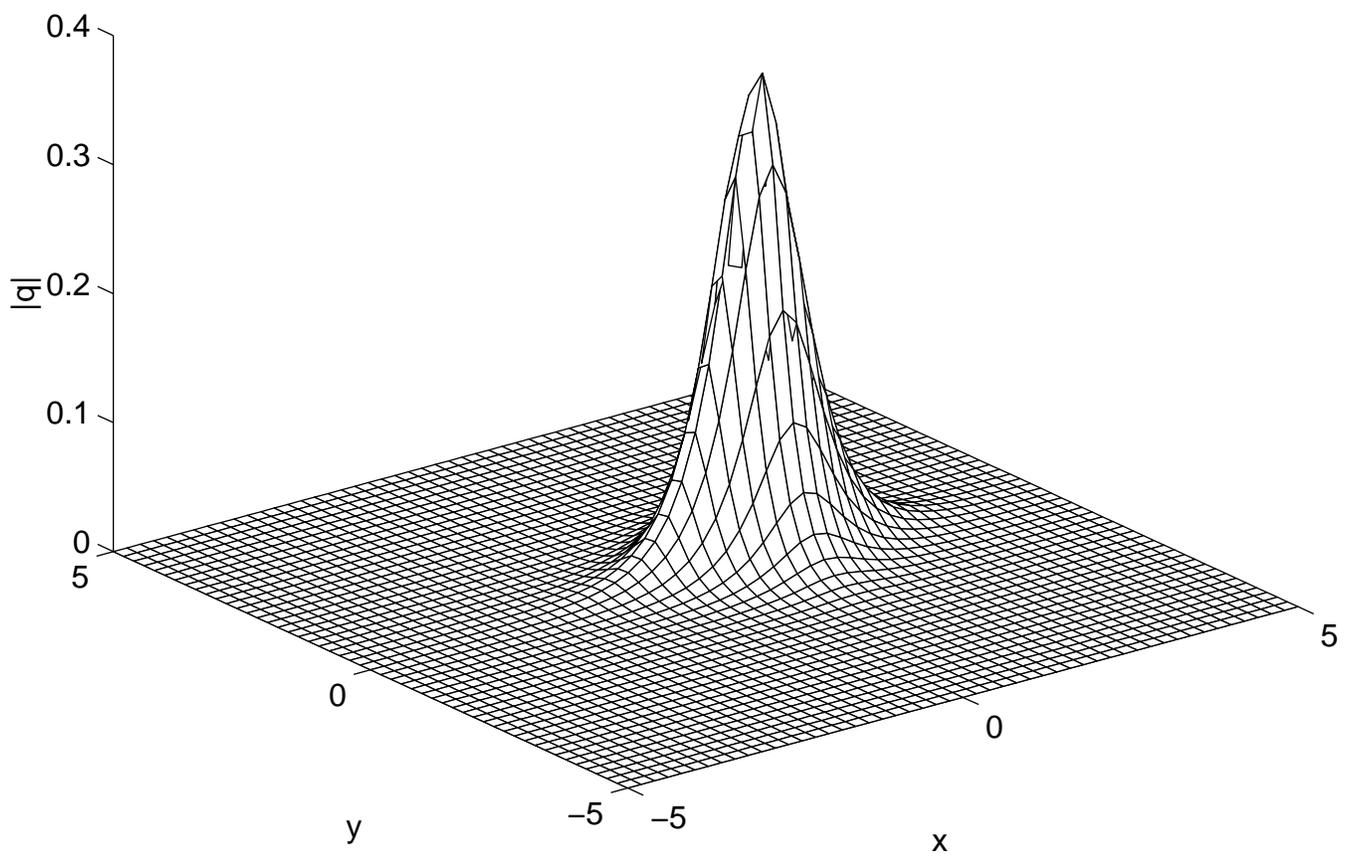
\end{document}